\documentclass[10pt,twocolumn,journal,twoside]{IEEEtran}
\usepackage{hyperref}
\usepackage{amssymb}

\newtheorem{corol}{Corollary}
\newtheorem{pro}{Proposition}

\begin{document}

\title{On minimal distance between $q$-ary bent functions }

\author{
\IEEEauthorblockN{ Vladimir N. Potapov}

\IEEEauthorblockA{Sobolev
Institute of Mathematics
\\ Novosibirsk, Russia \\ vpotapov@math.nsc.ru\\} }

\maketitle

\begin{abstract}
 The minimal Hamming distance  between distinct
$p$-ary bent functions of $2n$ variables is proved to be $p^n$ for
any prime $p$. It is shown that the number of $p$-ary bent functions
at the distance  $p^n$ from the quadratic bent function is equal to
$p^n(p^{n-1}+1)\cdots(p+1)(p-1)$ as $p>2$.


\end{abstract}

\section{Introduction}

 Bent functions is well known as Boolean functions with extremal nonlinear
 properties. Boolean bent function
 are intensively studied at present as they have numerous applications
 in cryptography, coding theory, and other areas.
  $q$-Ary  generalizations
 of bent functions are an interesting mathematical subject as well (see \cite{Tok2}).
 In this paper we consider Hamming distances between different bent
 functions and properties of bent functions at minimal distance from each
 other.
The Hamming distance between two discrete functions is the number of
arguments where these functions are differ. In other words, the
Hamming distance between two  functions $f$ and $g$ is the
cardinality of the support $\{x\in G \ |\ f(x)\neq g(x)\}$ of their
difference.

\section{Fourier transform on finite abelian groups}

Let  $G$ be a finite abelian group. Consider  a vector space $V(G)$
consisting of functions $f:G\rightarrow \mathbb{C}$ with inner
product
$$(f,g)=\sum\limits_{x\in G}f(x)\overline{g(x)}.$$
A function $f:G\rightarrow \mathbb{C}\backslash\{0\}$ mapping the
group to the non-zero complex numbers is called a character of $G$
if it is a group homomorphism from $G$ to $\mathbb{C}$, i.e.
$\phi(x+y)=\phi(x)\phi(y)$ for each $x,y\in G$. The set of
characters of an abelian group is an orthogonal basis of $V(G)$. If
$G=Z^n_q$  then we can define characters of  $G$ by equation
$\phi_z(x)=\xi^{\langle x,z\rangle}$, where $\xi=e^{2\pi i/q}$ and
$\langle x,y\rangle=x_1y_1+\dots+x_ny_n\, {\rm mod}\, q$ for each
$z\in Z^n_q$.  We may define the Fourier transform of a $f\in V(G)$
by the formula $\widehat{f}(z)=(f,\phi_z)/|G|^{1/2}$, i.e.,
$\widehat{f}(z)$ is the coefficients of the expansion of $f$ in the
basis of characters. Parseval's identity $(f,f)=\|f\|^2=
\|\widehat{f}\|^2$ and the Fourier inversion formula
$\widehat{(\widehat{f(x)})}= f(-x)$ hold. A proof of the following
equation there can be found in  \cite{Tao}.

\begin{pro}(uncertainty principle)\label{sbent1} For every $f\in V(G)$
the following inequality is true:
\begin{equation}\label{ebent1}
 |{\rm
supp}(f)||{\rm supp}(\widehat{f})|\geq |G|.
\end{equation}
 \end{pro}

 If $H $ is any subgroup of $G$, and we set $f$ to be the
characteristic function of $H$, then it is easy to see that $|{\rm
supp}(f)|=|H|$ and $|{\rm supp}(\widehat{f})| = |G|/|H|$, so
(\ref{ebent1}) is sharp. One can show that up to the symmetries of
the Fourier transform (translation, modulation, and homogeneity)
this is the only way in which (\ref{ebent1}) can be obeyed with
equality.

 If $p$ is prime then
 $Z^n_p$ can be considered as $n$-dimensional vector space
 over $GF(p)$.

\begin{corol}\label{cbent1}
Let  $p$ be a prime number. An equation $|{\rm supp}(f)||{\rm
supp}(\widehat{f})|=p^n$ holds iff $f=c\phi_z\chi^\Gamma$, where
$z\in G$, $c\in \mathbb{C}$ is a constant and $\chi^\Gamma$ is the
characteristic function of an affine space  $\Gamma$ in  $Z^n_p$.
\end{corol}

The following equation  can be found in   \cite{Sar} and \cite{Tsf}.
\begin{pro}\label{sbent2}
  If  $p$ is a prime number and  $\Gamma$ is a linear subspace in
  $Z^n_p$,
  then it holds
$$\sum\limits_{y\in \Gamma}\widehat{f}(y)=p^{{\rm dim}(\Gamma) -n/2}\sum\limits_{x\in
\Gamma^\bot}{f}(x).$$
\end{pro}

Define the convolution of  $f\in V(G)$ and $g\in V(G)$  by equation
 $f*g(z)=\sum\limits_{x\in G}f(x)g(z-x)$. It is well known that

\begin{equation}\label{ebent2}
\widehat{f*g}=|G|^{1/2} \widehat{f}\cdot \widehat{g}.
\end{equation}

We may define the Walsh--Hadamard transform of function
$g:Z^n_q\rightarrow Z_q$  by the formula
$W_g(z)=\widehat{\xi^g}(z)$.

\section{Bent functions}

A function $f:Z^n_q\rightarrow Z_q$ is called a $q$-ary bent
function iff  $|W_f(y)|=1$ for each $y\in Z^n_q$ or
$\widehat{\xi^f}\cdot\overline{\widehat{\xi^f}}=I$, where $I$ is
equal to
 $1$ everywhere (see \cite{KSW}, \cite{Tok1}). By using
(\ref{ebent2}) we can obtain that the definition of bent function is
equivalent to the equation
${\xi^f}*\overline{{\xi^f}}=|G|\chi^{\{0\}}.$ Then the matrix
$B=(b_{z,y})$, where $ b_{z,y}=\xi^{f(z+y)}$, is a generalized
Hadamard matrix.

 A bent function  $b$ is called regular iff there exists a function
$b':Z^n_q\rightarrow Z_q$ such that  $\xi^{b'}=\widehat{\xi^{b}}$.
Then $b'$ is a bent function as well.  If  $q$ is a prime power and
$n$ is even, then each bent function is regular. We assume below
that $p$ is a prime number and $n$ is even.

\begin{pro}\label{sbent3}
1) $\sum\limits_{j=0}^{q-1}\xi^{kj}=0$ as $k\neq 0\, {\rm mod}\, q$;

2) if $q$ is a prime number then $\xi$ is not a root of rational
polynomial function of degree less than $q-1$.

\end{pro}

\begin{corol}\label{cbent2} For any two $p$-ary bent functions $b$
and $b'$, it holds $|{\rm supp}(\xi^b-\xi^{b'})|=|{\rm
supp}(\widehat{\xi^b}-\widehat{\xi^{b'}})|$.
\end{corol}

It is sufficient to show that there are   $\frac{p-1}{2}$ different
numbers of type $|\xi^i-\xi^j|^2$, ($i\neq j$), and these numbers
are independent over $\mathbb{Q}$. Then from Proposition
\ref{sbent1} and Corollary \ref{cbent1} we obtain

\begin{corol}\label{cbent3}
The Hamming distance between two bent function on  $Z^n_p$ is not
less than $p^{n/2}$. If it is equal to $p^{n/2}$, then the
difference between these functions is equal to $c\chi^\Gamma$, where
$c\in Z_p$ and $\Gamma$ is an $n/2$-dimensional affine subspace.
\end{corol}

From Proposition \ref{sbent2} one can assume that the following
statements is true.

\begin{corol}\label{cbent41}
 If a bent function $b:Z^n_p\rightarrow Z_p$ is an affine function on
 an affine subspace
 $\Gamma$, then ${\rm dim}\Gamma \leq n/2$.
\end{corol}

\begin{corol}\label{cbent4}
If a bent function $b:Z^n_p\rightarrow Z_p$ is an affine function on
 an $n/2$-dimensional affine subspace, then there exist  $p-1$ bent function which differ from
  $b$ only on this subspace.
\end{corol}

Corollaries  \ref{cbent3} -- \ref{cbent4} was proved in \cite{Car1}
in the case of $p=2$. In \cite{Pot12'} it was found   spectrum of
potential small distances (less than the doubled minimum distance)
between two bent functions in the binary case.

\section{Quadratic forms}

A quadratic form $Q:(GF(q))^n\rightarrow GF(q)$ is called
non-degenerate iff   $\{x\in (GF(q))^n : \forall y\in (GF(q))^n\,
Q(y+x)=Q(y)\}=\{0\}$. A linear subspace $U$ in $(GF(q))^n$ is called
totally isotropic iff $Q(U)=0$. The maximal dimension of a totally
isotropic subspace is often called the Witt index of the form.
 If
$n=2d$, then the maximal Witt index of a non-degenerate forms of
degree $n$ is equal to $d$. All non-degenerate forms with the
maximal Witt index are equivalent. One of such quadratic forms is
determined by the equation $Q_0(v_1,\dots,v_d,u_1,\dots u_d)=
v_1u_1+\dots+v_du_d$. It is well known that $Q_0$ is a bent function
from Maiorana--McFarland  class (see \cite{Tok1}). The following
proposition is proved, for example, in  \cite{BCN} (p.274, Lemma
9.4.1)).

\begin{pro}\label{sbent4}
The number of the totally isotropic subspaces of $Q_0$  is equal to
$\prod\limits_{i=1}^{d}(q^{d-i}+1)$.

\end{pro}

It is easy to see that if $Q_0$ is an affine function  on a some
affine subspace, then it is an affine function on every coset.
Moreover, if $Q_0$ is an affine function  on a linear subspace of
dimension $d$, then this subspace is isotropic (here we assume that
$q>2$). Thus $Q_0$ is an affine function on all cosets of totally
isotropic subspaces and  it is not affine  on other linear subspaces
of dimension $d$.

From Proposition \ref{sbent4} and Corollary \ref{cbent4} we can
conclude that

\begin{corol}\label{cbent5}
If $p$ is a prime number and $p>2$, then there are
$p^d(p^{d-1}+1)\cdots(p+1)(p-1)$ $p$-ary bent functions  at the
distance  $p^{d}$ from  $Q_0$.
\end{corol}

In the binary case the analogous statement was proved in
\cite{Kol1}. In \cite{Kol2} was established that this bound  for the
number of bent functions at the minimal distance is reached only for
quadratic bent functions. It is natural to assume that this property
of $p$-ary quadratic bent functions is true for any prime $ p>2$.

\end{document}